\documentclass[letterpaper,12pt]{article}
\usepackage{graphicx}
\usepackage{amsthm, amsmath, amssymb}
\usepackage{mathrsfs,latexsym}
\usepackage{verbatim}
\newtheorem{thm}{Theorem}[section]
\newtheorem{cor}[thm]{Corollary}
\newtheorem{lemma}[thm]{Lemma}
\theoremstyle{remark}
\theoremstyle{definition}
\newtheorem{defn}[thm]{Definition}

\numberwithin{equation}{section}
\theoremstyle{definition}
\newtheorem{exmp}[thm]{Example}

\begin{document}
%
\title{Representing the Smallest Ideal and it's closure in the semigroup $ (\beta N ,.)$ as equivalence classes under the extension of divisibility relations}
\author{Salahddeen Khalifa\\Department of Mathematics and Computer Sience\\University of Missouri- St. Louis\\St.Louis, MO 63121 U.S.A\\e-mail: skkf57@mail.umsl.edu}
%
%
\maketitle
\begin{abstract} 
In this paper we will prove that all the elements in the smallest ideal $ K (\beta N)$ in the semigroup of the Stone Cech compactification $ (\beta N , .)$ of the discrete semigroup of natural numbers $ N $ under multiplication constitute a single equivalence class under the relation $ =_ m ( x =_m y $ if $ x |_m y $ and $ y |_m x)$. And all the elements in the closure of the smallest ideal $ Cl ( K (\beta N)) $ constitute a single equivalence class under the relation $ =_\sim ( x =_\sim y $ if $ x \tilde | y $ and $ y \tilde | x)$
\end{abstract}
\section{Introduction}
 we consider the discrete space $ N$ of natural numbers and it's Stone-Cech compactification $ \beta N = \{ x : x $   is ultrafilter on  $ N \}$ with topology that is introduced by a base $ \mathcal{B} = \{ \bar A : A \subseteq N \}$ where $ \bar A = \{ x \in \beta N : A \in x \}$ is a clopen subset of $ \beta N $ . 
Then according to ([1] Chapter 4 (4.1)) we can extend the usual multiplication (.) which is defined on the semigroup $ N $ to an operation that is defined on $ \beta N $ as follows:  
$ A \in x . y $ iff $ \{ n \in N : A / n \in y \} \in x $ where $ A / n = \{ m \in N : m n \in A \} $. $ \beta N $ contains the set of natural numbers $ N$ as principal ultrafilters which are defined on $ N$, and the set of non-principal ultrafilters are denoted by $ N^* = \beta N-N$. $| \beta N|=2^c$ where $ c= |R|$. \newline 
$ \mu = \{ A \subseteq N : A = A \uparrow \}$ is called the collection of \textbf{upward closed} subsets of $ N$ where $ A \uparrow = \{ n \in N : \exists a \in A , a|n \}$ and $ \nu = \{ A \subseteq N : A = A \downarrow \}$ is called the collection of \textbf{downward closed} subsets of $ N$ where $ A \downarrow = \{ n \in N : \exists a \in A , n | a \}$. $ \mu $ and $ \nu$ have the finite intersection property and if $ A \in \mu $ then $ N-A \in \nu $. 
\begin{lemma} 
(1) ([1] Theorem 4.15 ) Let $ x,y \in\beta N$ and let $ A\subseteq N. $ Then $ A \in x.y$ if and only if there exists $ B \in x $ and an indexed family $ (C_n)_{n \in B} $ in $ y $ such that $ \bigcup\limits_{n \in B} n C_n \subseteq A. $ \newline
(2)Let $ x,y \in\beta N$ and let $ A, B\subseteq N$. Then \newline
(b) If $ A \in x , B \in y $ , then $ A B \in xy $ \newline 
(c) $ \bar A . \bar B \subseteq \overline {A B} $ for any $ A,B\subseteq N $ 
\end{lemma}
\section{Ideals} 
A subset $ L$ of $ \beta N $ is a left ideal if $ \beta N . L \subseteq L $ , and $L$ is a minimal left ideal if for any left ideal $L_1 $ we have $ L \subset L_1 $. 
A subset $ R$ of $ \beta N $ is a right ideal if $ R.\beta N \subseteq R$ , and $ R$ is a minimal right ideal if for any right ideal $R_1$ we have $ R \subset R_1 $ , A subset $J$ of $ \beta N $ is an ideal if $J$ left (right) ideal , and $J$ is a minimal ideal if it is a minimal left (right) ideal. \newline 
Each set of the form $ \beta N x = \{ p x : p \in \beta N \}$ , where $ x \in \beta N $ is called principal left ideal generated by $ x$ , each set of the form $ x \beta N = \{ x p : p \in \beta N \}$ where $ x \in \beta N$ is called principal right ideal generated by $x$ , and each set of the form $ \beta N x \beta N = \{ p x q : p , q \in \beta N \}$ where $ x \in \beta N $ is called principal ideal generated by $x$ . Each minimal left ideal $L$ is principal $ : L = \beta N x$ for every $ x \in L$ , each minimal right ideal $R$ is principal $ : R = x \beta N $ for every $ x \in R $ , and each minimal ideal $J$ is principal: $ J = \beta N x \beta N $ for every $ x \in J$ 
\begin{thm}
	(a) ([1] Theorem 2.7 (d)) If $L$ is a left ideal and $R$ is a right ideal then $ L \cap R \neq \phi $ \newline 
	(b) ([1] Corollary 2.6, Theorem 2.7(a)) Every left (right) ideal in $ \beta N$ contains a minimal left (right) ideal. \newline 
\end{thm}
\begin{thm}
([1] Theorem 6.30, Corollary 6.41)There are $ 2^{c}$ disjoint minimal left (right) ideals in $\beta N $  and each of them contains $ 2^{c}$ elements. 
\end{thm} 
\begin{thm} 
(a) ([1], Theorem 1.51) $ \beta N $ has the smallest ideal denoted by $ K (\beta N)$ that is contained in all other ideals in $ \beta N $ . \newline 
(b)([1], Theorem 1.64) 
\begin {align*}
 K (\beta N) & = \cup \{ L : L  \text {is a minimal left ideal of}  \beta N \} \\
 & =  \cup \{ R : R \text{is a minimal right ideal of}  \beta N \} 
 \end {align*}
\end{thm}
\section{Divisibility relations on $ \beta N $ } 
We can extend the usual divisibility relation $|$ that is defined on $N$ to a divisibility relation on $ \beta N $ in many ways. 
If $ x , y \in \beta N $ then the most common extension divisibility relations of $|$ are the following $ :y$ is \textbf{left-divisible} by $ x , x |_l y $ if there is $ z \in \beta N $ such that $ y = z x$ . $ y $ is \textbf{right-divisible} by $ x , x |_r y , $ if there is $ z \in \beta N $ such that $ y = x z$ . $y$ is \textbf{mid-divisible} by $ x , x |_m y $ , if there are $ z , w \in \beta N $ such that $ y = z x w$ . $ y $ is \textbf{$ \tilde | $ -divisible} by $ x , x \tilde | y $ if for all $ A \subseteq N , A \in x $ implies $ | A = \{ n \in N : \exists a \in A , a |n \} \in y $. When $ x = n \in N $ , we write $ n | y , (y = n z , z \in \beta N)$ . 
\begin{lemma} 
(a)([3]) $ |_l , |_r \subset |_m \subset \tilde |$ \newline 
(b)([3]) All the relations $ |_l , |_r , |_m $ and $ \tilde | $ are preorders (reflexive and transitive), but none of them is symmetric and antisymmetric \newline Since $ |_l , |_r , |_m $ and $ \tilde | $ are not antisymmetric then we can introduce another relations on $ \beta N $ that will be symmetric such as $ : x =_l y $ if $ x |_l y $ and $ y |_l x , x =_r y $ if $ x |_r y $ and $ y |_r x , x =_m y $ if $ x |_m y $ and $ y |_m x$ and $ x =_\sim y $ if $ x \tilde | y $ and $ y \tilde | x $ 
\end{lemma}
\begin{lemma}
([4] Lemma 1.4)The following conditions are equivalent: \newline
(a) $ x\tilde | y $ \newline
(b) $ x\cap \mu\subset\ y\cap\mu $\newline
(c) $ y\cap\nu\subset\ x \cap\nu $
\end{lemma} 
\begin{lemma} 
(a) $ =_l , =_r , =_m $ and $ \tilde |$ are equivalence relations. \newline 
(b)([4] Lemma 1.5 (b)) For each $ x \in \beta N $ the sets $ x \uparrow_\sim = \{ y \in \beta N : x \tilde | y \}$ , $ x \downarrow_\sim = \{ y \in \beta N : y \tilde | x \} $ and $ [x]_\sim $ are closed subsets of $ \beta N$ . 
\end{lemma} 
\section{Levels of ultrafilters} 
\begin{defn} 
 ([5] Definition 2.4) Let $P $ be the set of prime numbers in $N$,  and let \newline 
$ L_0 = \{ 1\} $ \newline 
$ L_1 = \{ a : a \in P \}$ \newline 
$ L_2 = \{ a_1  a_2 : a_1 , a_2 \in P \} $ \newline 
. \newline 
. \newline 
. \newline 
$ L_n = \{ a_1 a_2 ......... a_n : a_1 , a_2 , ........, a_n \in P \}$. 
Then, the ultrafilters $x$ is called on " finite level" if it is in exactly on one of the following sets $ : \bar L_i , i = 0 , 1 , 2, .......$ where $ \bar L_i = \{ y \in \beta N : L_i \in y \}$ 
\end{defn} 
\begin{lemma}
	 
([2] Lemma 3.1(c)) $ \bigcup\limits^{\infty}_{i=0} \bar L_i \neq \overline {\bigcup\limits^{\infty}_{i=0}} L_i $ \newline
 By the lemma(4.2) there are ultrafilters that are not on finite levels, these ultraafilters are given by the following definition 
\end{lemma} 
\begin{defn} 
	( [2] Definition 3.2 ) The ultrafilters that belong to \newline $ \mathcal {I} = \bigcap\limits^{\infty}_{i=0} \overline {N-L_i} $ are called the "ultrafilters that are not on finite levels". 
\end{defn}
\begin{thm} 
([2] Theorem 3.7)\newline
(a) If $ x , y \in \beta N $ and $ x , y \in \mathcal {I} $ . Then $ x y \in\mathcal {I}  $ and $ y x \in \mathcal  {I} $ \newline 
(b) If $ x , y \in \beta N $ and $ x \in\mathcal {I}  , y \notin \mathcal {I}  $ . Then $ x y \in \mathcal {I} $ and $ y x \in \mathcal {I}  $ 
\end{thm} 
\section{Smallest ideal $ K (\beta N)$ in the semigroup $ (\beta N , .)$ }
\begin{lemma}
	(a) $ \mathcal {I}$  is a left ideal in the semigroup $ (\beta N,.)$ \newline 
	(b) $ \mathcal {I}$  is a right ideal in the semigroup $ (\beta N,.)$ \newline 
	(c) $ \mathcal {I}$ is an ideal in the semigroup $ (\beta N,.)$ \newline  
\end{lemma} 
\begin{proof}
	By (Theorem 4.4)  
\end{proof} 
\begin{lemma}
	(a) $ \mathcal {I}$ is not a minimal left ideal in $ (\beta N,.)$ \newline 
	(b) $\mathcal {I}$ is not a minimal right ideal in $ (\beta N,.)$ \newline 
	(c) $ \mathcal {I}$ is not a minimal ideal in $ (\beta N ,.)$  
\end{lemma} 
\begin{proof} 
	(a) Let $ L = \{ 2 x : x \in \mathcal {I} \}$  , then $ L$ is a left ideal and $ L \subset \mathcal {I} , L \neq \mathcal {I}$  [ by ( [1], Lemma 6.28 ) $ 2 x \neq n x $ for any $ n \in N , n \neq 2 $ ] . Thus $ \mathcal {I}$ is not a minimal left ideal. \newline 
	(b) Similar to (a) \newline 
	(c) By (a) and (b) 
\end{proof} 
\begin{cor} 
	The smallest ideal $ K (\beta N)$ of the semigroup $ ( \beta N,.)$ exists and is a subset of $ \mathcal {I} $ i.e $  K (\beta N ) \subset \mathcal {I} $ 
\end{cor} 
\begin{proof} 
	By (Theorem 2.2) $ \beta N $ has $ 2 ^c $ minimal left (right) ideals and by ( [1] Lemma 1-45) for any minimal left ideal $ L$ and minimal right ideal $R$ we have $ L \subseteq \mathcal{I} $ and $ R \subseteq \mathcal {I} $ , and by (Lemma 5.2) since $\mathcal {I} $ is not a minimal ideal, then  $ L \subset \mathcal {I} $ and $ R \subset \mathcal {I} $ . Then, by (Theorem 2.3 (a) ) $ K (\beta N) $  exists and 
	\begin{align*}
	  K ( \beta N) & = \cup \{ L : L \text { is a minimal left ideal of }  \beta N \}\\ 
	& = \cup \{ R : R \text { is a minimal right ideal of }  \beta N \} 
	\end{align*} 
	Thus $ K ( \beta N) \subset \mathcal {I}$ 
\end{proof} 
\begin{lemma}
	(a) If $ L $ is a minimal left ideal and $x\in L$. Then $ x |_l y $ if and only if $ y \in L $ \newline 
	(b) If $ R $ is a minimal right ideal and $ x\in R$, Then $ x |_r y $ if and only if $ y \in R $ \newline 
	(c) If $ J $ is a minimal ideal and $x\in J$. then $ x |_m y $ if and only if $ y \in J $  
\end{lemma}
\begin{proof}
	(a) $ ( \Rightarrow ) $ Let $ L$ be a minimal left ideal. Then $ L = \beta N x $  for any $ x \in L $. Let $ x |_l y $ , then $ y = z x $ for some $ z \in \beta N $ . Thus $ y \in \beta N x $ so $ y \in L \newline (\Leftarrow) $ If $ y \in L $ , then $ y \in \beta N x $ for any $ x \in L , $ so $ y = z x $ for some $ z \in \beta N $ Thus $ x |_l y $ \newline 
	(b) $ ( \Rightarrow ) $ Let $ R$ be a minimal right ideal. Then $ R = x \beta N $ for any $ x \in R $. Let $ x |_r y $ , then $ y = x z $ for some $ z \in \beta N. $ Thus $ y \in x \beta N $ , so $ y \in R \newline (\Leftarrow) $ If $ y \in R $ then $ y \in x \beta N $ for any $ x \in R $ , so $ y = x z$ for some $ z \in \beta N$ . Thus $ x |_r y $ \newline 
	(c) $ ( \Rightarrow ) $ Let $ J$ be a minimal ideal. Then $ J = \beta N x \beta N $ for any $ x \in J $. Let $ x |_m y $ , then $ y = z x w $ for some $ z ,w \in \beta N $ . Thus $ y \in \beta N x \beta N $ so $ y \in J \newline (\Leftarrow) $ If $ y \in J $ , then $ y \in \beta N x \beta N $ for any $ x \in J $ , so $ y = z x w $ for some $ z , w \in \beta N $ . Thus $ x |_ m y $
\end{proof}
\begin{thm} 
	Let $ (\beta N , =_l) $ where $ x =_l y $ if $ x |_l y $ and $ y |_l x$ . Then \newline 
	(a) If $ x \in K (\beta N) $ , then $ [x]_l = L , L $ is a minimal left ideal, $ x \in L $ , and $ [x]_l $ has $ 2^{c} $ elements. \newline 
	(b) $ K (\beta N) $ is partitioned by $ 2^{c}$ equivalence classes $[x]_l$.   
\end{thm}
\begin{proof} 
	(a) Let $ x \in K (\beta N)$ , then by (Theorem 2.2) there exists only one minimal left ideal $ L $ such that $ x \in L $. Since $ L$ is a minimal left ideal , then $ L$ is principal , so $ L = \beta N x$ . Also for any $ y \in L $ we have $ L = \beta N y $ . Then $ y \in \beta N x $ , and there exist $z\in\beta N$ such that $ y=zx$.Thus $x|_l y$. Also $x\in\beta Ny$, and there exist $ w \in \beta N $ such that $ x = w y $ . Thus $ y |_l x $ . Therefore $ x =_l y $ for any $ y \in L $ . If $ p \in \beta N , p \notin L = \beta N x $ , by (Lemma 5.4(a)) we have $ x \neq_l p $ . So $ [x]_l = L $ . Since $ L$ has $ 2^{c} $ elements , then $ [x]_l $ has $ 2^{c} $ elements. \newline 
	(b) Since $ K (\beta N)$ has $ 2^{c}$ disjoint minimal left ideals, and by (a) each of them is $ [x]_l $ equivalence class. Then $ K (\beta N)$ is partitioned by $2^c$ equivalence classes $[x]_l$ . 
\end{proof} 
\begin{thm} 
	Let $ (\beta N , =_r)$ such that $ x =_r y $ if $ x |_r y $ and $ y |_r x$ . Then \newline 
	(a) If $ x \in K (\beta N)$ , then $ [x]_r = R , R $ is a minimal right ideal $ , x \in R $ , and $ [x]_r $ has $ 2^{c}$ elements . \newline 
	(b) $ K (\beta N)$ is partitioned by $ 2^{c} $ equivalence classes$[x]_R$. 
\end{thm} 
\begin{proof} 
	Similar to (Theorem 5.5) 
\end{proof}
\begin{thm}
	Let $ (\beta N , =_m) $ where $ x =_m y $ if $ x |_m y $ and $ y |_m x $ . Then \newline 
	(a) If $ x \in K (\beta N) $ , then $ [x]_m = K (\beta N)$ , and $ [x]_m $ has $ 2^{c}$ elements . \newline 
	(b) $ =_m $ has one equivalence class for any $ x \in K (\beta N)$ .
\end{thm} 
\begin{proof}
	(a) By (Theorem 2.1 (a) ) for any minimal left ideal $ L$ and any minimal right ideal $ R$ we have $ L \cap R \neq \phi $ . So there exist $ x $ such $ x \in L $ and $ x \in R $. If $ x \in L , $ then for any $ y \in L$ we have $ y \in \beta N x$ , and $ x \in \beta N y $ , so $ x |_l y $ and $ y |_l x , x |_m y $ and $ y |_m x$ . Thus $ x =_m y $ for any $ y \in L $ . If $ x \in R $ , then for any $ z \in R $ we have $ z \in x \beta N$ and $ x \in z \beta N$ , so $ x |_R z$ and $ z |_R x , x |_m z$ and $ z |_m x$ . Thus $ x =_m z$ for any $ z \in R $ . If $ w \in \beta N, w \notin L$ and $ w \notin R,$ then $ w \notin \beta N x \beta N $ , and by (Lemma 5.4 (c)) we have $ w \neq_m x $ . So $ [x]_m = K (\beta N).$ Since $ K (\beta N)$ has $ 2^{c}$ elements , so $ [x]_m $ has $ 2^{c}$ elements . \newline 
	(b) By (a) we have one equivalence class such that $ | [x]_m | = 2^{c}$ for any $ x \in K (\beta N).$  
\end{proof}
\section{Closure of $ K (\beta N) $ in the semigroup $ (\beta N, .)$}
\begin{defn} 
	Consider the semigroup $ (N, .)$. Then \newline 
	(a)([1] Definition 4.45) A set $ A \subseteq N $ is "thick" if the family $ \{ A / n : n \in N \}$ has the finite intersection property. \newline 
	(b)([1] Definition 4.38(a)) A set $ A \subseteq N $ is "syndetic" if there is a finite subset $ E $ of $ N$ $ ( E \in \mathcal{P}_f (N))$ such that $ \bigcup\limits_{n \in E} {A /n} = N $. \newline 
	(c)([1]Definition 4.38(b)) A set $ A \subseteq N $ is piecewise syndetic if there is a finite subset $ E$ of $ N$ $ ( E \in \mathcal{P}_f (N))$ such that the set $ \bigcup\limits_{n\in E} A / n$ is thick in $ N$ i.e $ \{ ( \bigcup\limits_{n \in E} A / n ) / m : m \in N \}$ has the finite intersection property.  
\end{defn}
\begin{thm}
	Let $ A \subseteq N $ . Then \newline 
	(a)([1] Theorem 4.48(a)) A is thick if and only if there is a left ideal of $ \beta N$ contained in $ \bar A $ . \newline 
	(b)([1] Theorem 4.48(b)) A is syndetic if and only if for every left ideal $ L $ of $ \beta N , L \cap \bar A \neq \phi $ . \newline 
	(c)([1] Theorem 4.40) A is piecewise syndetic if and only if $ \bar A \cap K (\beta N) \neq \phi$.  
\end{thm} 
\begin{lemma}
	([1] 4.46) \newline
	(a) Every thick set is piecewise syndetic . \newline 
	(b) Every syndetic is piecewise syndetic. 
\end{lemma}
\begin{exmp} 
	(1)$ N-L_i $  is thick for all $ i = 0 , 1 , .......$ since for any minimal left ideal $ L$ of $ \beta N $ we have $ L \subseteq \overline {N-L_i}$ . Thus by (Theorem 6.2(a)) $ N-L_i$ is thick. \newline 
	(2)$ N-L_i $ is syndetic for all $ i = 0 , 1, ........$ Since $ \overline {N-L_i} \cap L \neq \phi$ for every left ideal $ L$ of $\beta N $ . Then by (Theorem 6.2(b)) $ N-L_i $ is syndetic. \newline 
	(3) $ N-L_i $ is piecewise syndetic.   
\end{exmp}  
\begin{lemma} 
	(a) $ A\subseteq N $ is syndetic if and only if $ N-A$ is not thick \newline 
	(b) $ A\subseteq N$ is thick if and only if $ N-A$ is not syndetic
\end{lemma} 
\begin{lemma}
	Let $ \mu $ be the set of all upward closed subsets of $N$ and $ \nu$ be the set of all downward closed subsets of $N$. Then \newline 
	(a) Any $ A \in \mu $ is thick \newline 
	(b) Any $ A \in \nu$ is not thick \newline 
	(c) Any $ A \in \mu$ is syndetic \newline 
	(d) Any $ A \in \nu $ is not syndetic  
\end{lemma}
\begin{proof} 
	(a) Let $ A \in \mu $.  \newline 
	$ A/2 = \{ m \in N : 2 m \in A \}, $ and  
	$ A /2 \uparrow = \{ n \in N : \exists a \in A /2 , a|n \} $. \newline 
	Let $ n \in A/2 \uparrow $ , so $ n = r a $ for some $ r\in N $ and $ a \in A / 2 $ , so $ 2 a \in A $. Since $ A \in \mu $ then $ r (2 a) \in A $ $ (2 a\mid r(2a)) $ , so $ 2 (r a) \in A$ , and $ r a \in A / 2$ . Thus $ n \in A/2 , $ and since $ A /2 \subseteq A/2 \uparrow $ . So $ A /2 = A /2 \uparrow $ i.e $ A / 2\in \mu$ \newline In general. Let $A\in\mu$. \newline 
	$ A / n = \{ m \in N : m n \in A \}, $ and 
	$ A/ n \uparrow = \{ r \in N : \exists a \in A / n , a |r \}$. \newline Let $ r \in A / n \uparrow $ , so $ r = ta $ for some $ t\in N $ and $ a \in A / n$ , so $ n a \in A$ . Since $ A$ is upward closed, then $ t ( n a) \in A \left( n a | t ( n a) \right) $ , so $ n (t a) \in A $  , and $ t a \in A / n$ . Thus $ r \in A /n $, and since $ A /n \subseteq A /n \uparrow $ , then $ A / n = A /n \uparrow $ , so $ A / n \in \mu $ . Therefore $ \{ A /n : n \in N \} \subseteq \mu$ . And since $\mu$ has $ F.I.P$ , so $ \{ A /n : n \in N \}$ has $ F.I.P$ , i.e $ A$ is thick \newline 
	(b) Let $ A \in \nu , A /m = \{ r \in N : r m \in A \}$. If $ m \notin A$ then $ A / m = {\{ r \in N : r m \in A  }\}= \phi $ , because if $ r m \in A$ , and since $ m | r m $ , then $ m \in A $ a contradiction . So, if $ A \in \nu $ , then $ \{ A /n : n \in N \}$ does not have $ F.I.P$ . Therefore $A$ is not thick \newline 
	(c) Let $ A \in \mu $ , we have $ N-A \in \nu $ and by (b) $ N-A $ is not thick . Thus by (Lemma 6.5(a)) $ A$ is syndetic . \newline 
	(d) Let $ A \in \nu $ , then we have $ N-A \in \mu $ , and by (a) $ N-A$ is thick . Thus by (Lemma 6.5 (b)) $A$ is not syndetic .  	 
\end{proof}
\begin{defn} 
	The set  of ultrafilters which are generated by the set of upward closed subsets of $ N $ is denoted by $ \hat {\mu}$ and it is given by $ \hat {\mu} = \bigcap\limits_{A \in \mu} \bar A $ \newline i.e $  \hat {\mu} = \{ x \in \beta N : x\in\bar A, A\in\mu \} $. 
\end{defn}
By (Lemma 3.2) the set of ultrafilters which belong to $ \hat {\mu}$ are divided by all ultrafilters of $ \beta N $ under $ \tilde\mid$- divisibility
\begin{lemma} 
	(a) For any ultrafilter $ x \in \beta N $ such that $ x \in \hat {\mu}, L_i \notin x$ $( x \notin \overline L_i)$ for all $ i = 0 , 1,2,......$ \newline 
	(b) $ x \in \mathcal{I} $ for any $ x \in \hat {\mu}$ $( \hat {\mu} \subset\mathcal {I})$ \newline 
	(c) There are $ 2^{c}$ ultrafilters $ x $ such that $ x \in\mathcal {I} , x \notin \hat {\mu} $ 
\end{lemma} 
\begin{proof}
	(a) Let $ x \in \beta N , x \in \hat {\mu} $ . First we need to prove that any subset $ A \subseteq N, A \in x$ where $ x \in \hat {\mu}$ is syndetic. For any element $ A \in x$ we have three cases: $ A \in \mu $ or $ A$ is a finite intersection of elements of $ \mu$ or $A $ contains some elements from $ \mu$ . \newline 
	(1) If $ A \in \mu$ then by (Lemma6.6(c)) $ A$ is syndetic \newline 
	(2) If $ A$ is finite intersection of elements of $\mu$ , then $ A \in \mu $ and $ A$ is syndetic. \newline 
	(3) If $ B \subset A, B \in \mu $ , since $ B$ is syndetic , then by (Theorem 6.2 (b) ) $ L \cap \bar B \neq \phi $ for any left ideal $ L$ of $ \beta N$ , so $ L \cap \bar A \neq \phi $ . Thus $ A$ is syndetic .\newline We know, by (Example 6-4) $ L_i$ is not syndetic for any $ L_i , i = 0 , 1,2.......$ Therefore $ L_i \notin x$ $ (x \in \bar L_i)$ \newline 
	(b) Let $ x \in \beta N$ , and $ x \in \hat {\mu}$ by (a) $ L_i \notin x$ for all $ i =0, 1,.........,$ so $ N-L_i \in x , x \in \overline {N-L_i} .$ Therefore $ x \in \bigcap\limits^\infty_{i=0} \overline{N-L_i} $ , Thus $ x \in \mathcal {I} $ and $ \hat {\mu}\subset \mathcal {I}$ \newline 
	(c) Let $ A = \{ n_1, n_2 ,n_3, ........ \} $ such that $ n_i \in L_i , i \geq 1 , n_i $ are odd number,  so by ([6] Theorem 3.3 ) $ \bar A$ contains $ 2^{c}$ ultrafilters $ x$ and since there is no nonprincipal ultrafilters $ y $ such that $ y\in\bar A$ and $ y\in\bar L_i$ for all $ i= 1,......$ Thus $ x \in \mathcal {I} $ . And let $ B = \{ 2, 4,6,.....\}$ , then $ B \in \mu$ and $ B \notin x $ for any $ x \in \bar A$ $(A \cap B = \phi)$. Thus $ x \notin \hat {\mu} $ for any $ x \in \bar A $.
\end{proof}
\begin{lemma} 
	(a) $ \hat {\mu}$ is a left ideal in the semigroup $ (\beta N,.)$ \newline 
	(b) $ \hat {\mu}$ is a right ideal in the semigroup $ (\beta N,.)$ \newline 
	(c) $ \hat {\mu}$ is an ideal in the semigroup $ (\beta N,.)$ . \newline 
	(d) $ K (\beta N) \subset \hat {\mu}$ \newline 
	(e) $ |\hat {\mu}| = 2^{c}$ 
\end{lemma} 
\begin{proof} 
	(a)Let $ x \in \hat {\mu} $ and $ y \in \beta N$ . For any $ A \in \mu $ we have $ A \in x $ , and for any $ B \in y $ by ( Lemma 1.1 (2) ) we have $ y x \in \bar B \bar A \subseteq \overline {A B}$ , so $ A B \in y x $. Since $ A \in \mu $ , then for any $ b a \in B A$ we have $ b a \in A $ , so $ B A \subseteq A$ . Thus $ A \in y x $ for any $ A \in \mu $ . Therefore $ y x \in \hat {\mu}$ . \newline 
	(b) Similar to (a) \newline 
	(c) By (a) and (b) \newline 
	(d) By (c) \newline 
	(e) Since $ |K (\beta N)| = 2 ^{c}$ , then by (d) we have $ |\hat {\mu}| = 2^{c}$.  
\end{proof} 
\begin{lemma} 
	(a) Let $ (\beta N ,  =_\sim) $ where $ x =_\sim y $ if $ x \tilde | y$ and $ y \tilde | x$ . Then $ [x]_\sim = \hat {\mu} $ for any $ x \in \hat {\mu} $ \newline 
	(b) $ K (\beta N) \neq \hat {\mu}$ 
\end{lemma} 
\begin{proof} 
	(a) Let $ x \in \hat {\mu} $ , then for any $ y \in \hat {\mu} $ we have $ x \cap \mu \subseteq y $ and $ y \cap \mu \subseteq x$, so by (Lemma 3.2 ) we have $ x \tilde | y $ and $ y \tilde | x$ . So $ x =_\sim y $ for any $ y \in \hat {\mu} $ . If $ y\notin \hat {\mu}$, then $ x\cap\mu\not\subseteq y$ . So $ x\not\tilde|  y $. Therefore $ [x]_\sim = \hat {\mu} $ for any $ x \in \hat {\mu} $ \newline 
	(b) Since for any $ x \in k (\beta N), $ by (Theorem 5.7 (a)) we have $ [x]_m = K (\beta N) , $ and by (a) we have $ [x]_\sim = \hat {\mu} $ , and since $ |_m \subset \tilde |$ , then $ [x]_m \subset [x]_\sim$ . Thus $ K (\beta N) \neq \hat {\mu}$ 
\end{proof}
\begin{cor}
	$ \hat {\mu} $ is closed subset of $ \beta N$.  
\end{cor} 
\begin{proof} 
	By (Lemma 3.3 (b)) $ [x]_\sim $ is closed and compact for any $ x \in \beta N$ . Then by (Lemma 6.10) $ \hat {\mu} $ is closed subset of $ \beta N $. 
\end{proof}
\begin{lemma}
	Any subset $ A \subseteq N$ such that $ A \in x $ where $ x \in \hat {\mu} $ is piecewise syndetic  
\end{lemma}
\begin{proof} 
	Let $ A \subseteq N$ such that $ A \in x , x \in \hat {\mu} $ . Analogues to proof of ( Lemma 6.8(a)) we have $A$ is syndetic , and by (Lemma 6-3) $A$ is piecewise syndetic. 
\end{proof} 
\begin{thm}
	(a)  If $ x \in \hat {\mu} $ . Then $ x \in Cl ( K (\beta N))$. \newline 
	(b) $ Cl (K(\beta N)) = \hat {\mu}$. \newline 
	(c) $ =_\sim $ has one equivalence class for any $ x \in Cl ( K (\beta N))$.
\end{thm}
\begin{proof}
	Let $x\in \hat {\mu}$ and let U be any open subset of $\beta N$ such that $x\in U$. So there exist a basic open subset $ \bar A$ of $\beta N$ such that $x\in\bar A\subset U$,and since $A\in x$ then by (Lemma 6.12) A is piecewise syndetic.So by (Theorem 6-2)we have $U\cap\ K(\beta N)\neq\phi $.Thus $x\in Cl(K(\beta N))$. \newline 
	(b) By (Lemma 6.9 (d)) we have $ K(\beta N)\subseteq\hat{\mu}$,so $ Cl(K(\beta N))\subseteq\ CL(\hat{\mu})$ and by (Corollary 6.11) we have $ Cl(K(\beta N)) \subseteq\hat{\mu}$. Since by (a) $ \hat{\mu}\subseteq\ Cl(K(\beta N))$. Thus $ Cl((K\beta N))=\hat{\mu}$. \newline
	(c) By (b) and ( Lemma 6.10(a) )
\end{proof}

\end{document}